\documentclass[12pt,a4paper]{amsart}
\usepackage{amssymb,amsmath}

\headsep=1cm \numberwithin{equation}{section}
\newtheorem{theorem}{Theorem}[section]
\newtheorem{lemma}[theorem]{Lemma}
\newtheorem{proposition}[theorem]{Proposition}
\newtheorem{corollary}[theorem]{Corollary}

\theoremstyle{definition}
\newtheorem{definition}[theorem]{Definition}
\theoremstyle{remark}

\newtheorem{example}[theorem]{Example}

\newcommand{\Ass}{\operatorname{Ass}}

\newcommand{\Spec}{\operatorname{Spec}}

\newcommand{\V}{\operatorname{V}}
\newcommand{\Z}{\operatorname{Z}}
\newcommand{\E}{\operatorname{E}}

\newcommand{\Ext}{\operatorname{Ext}}

\newcommand{\Hom}{\operatorname{Hom}}

\newcommand{\Ann}{\operatorname{Ann}}

\newcommand{\Gdim}{\operatorname{Gdim}}

\newcommand{\lo}{\longrightarrow}
\newcommand{\fm}{\frak{m}}
\newcommand{\fp}{\frak{p}}
\newcommand{\fq}{\frak{q}}
\newcommand{\fa}{\frak{a}}

\newcommand{\fn}{\frak{n}}

\begin{document}
\author[Divaani-Aazar and Esmkhani]{Kamran Divaani-Aazar and Mohammad Ali
Esmkhani}
\title[Artinianness of local cohomology ... ]
{Artinianness of local cohomology modules of ZD-modules}

\address{K. Divaani-Aazar, Department of Mathematics, Az-Zahra University,
Vanak, Post Code 19834, Tehran, Iran-and-Research Institute for
Fundamental Sciences, Tabriz, Iran.} \email{kdivaani@ipm.ir}

\address{M.A. Esmkhani, Department of Mathematics, Shahid Beheshti University,
Tehran, Iran.}


\keywords{Local cohomology, Artinian modules, $ZD$-modules, Goldie
dimension.\\
This work has been supported by the Research Institute for
Fundamental Sciences, Tabriz, Iran.}

\begin{abstract}

This paper centers around Artinianness  of the local cohomology of
$ZD$-modules. Let $\fa$ be an ideal of a commutative Noetherian
ring $R$. The notion of $\fa$-relative Goldie dimension of an
$R$-module $M$, as a generalization of that of Goldie dimension is
presented. Let $M$ be a $ZD$-module such that $\fa$-relative
Goldie dimension of any quotient of $M$ is finite. It is shown
that if $\dim R/\fa=0$, then the local cohomology modules
$H^i_{\fa}(M)$ are Artinian. Also, it is proved that if $d=\dim M$
is finite, then $H^d_{\fa}(M)$ is Artinian, for any ideal $\fa$ of
$R$ . These results extend the previously known results concerning
Artinianness  of local cohomology of finitely generated modules.
\end{abstract}

\maketitle

\section{Introduction}

Throughout this paper, $R$ is a commutative Noetherian ring with
identity and all modules are assumed to be unitary. Let $\fa$ be
an ideal of $R$ and $M$ an $R$-module. The $\fa$-torsion submodule
$\cup_{n\in \mathbb{N}}(0:_M\fa^n)$ of $M$ is denoted by
$\Gamma_{\fa}(M)$. For each integer $i\geq 0$, the i-th local
cohomology functor $H^i_{\fa}(.)$ is defined as the i-th right
derived functor of $\fa$-torsion functor $\Gamma_{\fa}(.)$. Also,
it is known that for each  $i\geq 0$ there is a natural
isomorphism of $R$-modules
$$H^i_{\fa}(M)\cong\underset{n}{\varinjlim}\Ext^i_R(R/\fa^n,M).$$
We refer the reader to text book [{\bf 2}] for more details about
local cohomology.

It is known that the local cohomology of finitely generated
modules have many interesting properties. In particular, if
$(R,\fm)$ is a local ring and $M$ a finitely generated $R$-module,
then the local cohomology modules $H^i_{\fm}(M)$ are Artinian.
Also, in the same situation, it is known that for $d=\dim M$, the
$d$-local cohomology module of $M$ with respect to any ideal $\fa$
is Artinian. It will be a noticeable achievement, if we could
extend these results to local cohomology of a larger class of
modules. In this paper, we shall show that $ZD$-modules behave
very well in conjunction with Artinianness  of local cohomology
modules.

An $R$-module $M$ is said to be $ZD$-module (zero-divisor module)
if for any submodule $N$ of $M$, the set of zero divisors of $M/N$
is a union of finitely many prime ideals in $\Ass_R(M/N)$.
According to Example 2.2, the class of $ZD$-modules is much larger
than that of finitely generated modules. As the main result of
this paper, we prove that for a $ZD$-module $M$ the
following are equivalent:\\
i) $\Gamma_{\fa}(M/N)$ is Artinian for any submodule $N$ of
$M$.\\
ii)$H^i_{\fa}(M/N)$ is Artinian for any submodule $N$ of $M$ and
all $i\geq 0$.

We say that an $R$-module $M$ has finite $\fa$-relative Goldie
dimension if the Goldie dimension of the $\fa$-torsion submodule
of $M$ is finite. Clearly, $\fa$-relative Goldie dimension of any
finitely generated module is finite. Let $M$ be a $ZD$-module such
that $\fa$-relative Goldie dimension of any quotient of $M$ is
finite. By using the above mentioned result, we deduce that the
local cohomology modules $H^i_{\fa}(M)$ are
Artinian if either,\\
i) $\dim R/\fa=0$, or\\
ii) $d=\dim M$ is finite and $i=d$.

\section{$ZD$-modules and Goldie dimension}

Let for an $R$-module $M$, $\Z_R(M)$ denote the set of zero
divisors on $M$. Evans [{\bf 4}] calls a ring $R$ a $ZD$-ring
(zero-divisor ring) if for any ideal $\fa$ of $R$, $\Z_R(R/\fa)$
is a union of finitely many prime ideals. Next, we present the
following modification of the definition of $ZD$-modules in [{\bf
6}].
\begin{definition} An $R$-module $M$ is said to be $ZD$-module if
for every submodule $N$ of $M$, the set $\Z_R(M/N)$ is a union of
finitely many prime ideals in $\Ass_R(M/N)$.
\end{definition}

An $R$-module $M$ is said to be Laskerian if any submodule of $M$
is an intersection of a finite number of primary submodules.
Obviously, any Noetherian module is Laskerian. An $R$-module $M$
[{\bf 3}] is said to be weakly Laskerian if the set of associated
primes of any quotient module of $M$ is finite. Clearly, any
Laskerian module is weakly Laskerian and any weakly Laskerian
module is $ZD$-module. In the sequel, we provide a large variety
of examples of $ZD$-modules.

\begin{example} i) It is easy to see  that any module with finite support
is weakly Laskerian. In particular, any Artinian module is a $ZD$-module.
 Also, by using this fact we can provide examples of $ZD$-modules which
are neither finitely generated nor Artinian.\\
ii) Recall that a module $M$ is said to have finite Goldie
dimension if $M$ does not contain an infinite direct sum of
non-zero submodules, or equivalently, the injective envelope
$\E(M)$ of $M$ decomposes as a finite direct sum of indecomposable
injective submodules. Because for any $R$-module $C$, we have
$\Ass_RC=\Ass_R\E(C)$, it turns out that any module with finite
Goldie dimension has only finitely many associated prime ideals.
This yields that a module of which all quotients have
finite Goldie dimension is weakly Laskerian.\\
iii) Let $E$ be the minimal injective cogenerator of $R$ and $M$
an $R$-module. If for an $R$-module $M$ the natural map from $M$
to $\Hom_R(\Hom_R(M,E),E)$ is an isomorphism, then $M$ is said to
be Matlis reflexive. By [{\bf 1}, Theorem 12], an $R$-module $M$
is Matlis reflexive if and only if $M$ has a finitely generated
submodule $S$ such that $M/S$ is Artinian and $R/\Ann_RM$ is a
complete semi-local ring. Also, by [{\bf 5}, Corollary 1.2], any
quotient of an $R$-module $M$ has finite Goldie dimension if and
only if $M$ has a finitely generated submodule $S$ such that $M/S$
is Artinian. Thus, by
(ii) any Matlis reflexive module is a $ZD$-module.\\
iv) An $R$-module $M$ is said to be linearly compact if each
system of congruences $$x\equiv x_i (M_i),$$ indexed by a set $I$
and where the $M_i$ are submodules of $M$, has a solution $x$
whenever it has a solution for every finite subsystem. It is clear
that, every quotient of a linearly compact module is also linearly
compact. On the other hand a linearly compact module $M$ has
finite Goldie dimension (see e.g. [{\bf 9}, Chapter 1.3]). Thus,
if $M$ is a linearly compact module, then any quotient of $M$ has
finite Goldie dimension, and so, by (ii) $M$ is a $ZD$-module.
\end{example}

Next, we bring the following characterization of $ZD$-modules.

\begin{lemma} Let $M$ be an $R$-module. The following are equivalent:\\
i) $M$ is a $ZD$-module.\\
ii) For every submodule $N$ of $M$, the number of prime ideals
with the property being maximal in $\Ass_R(M/N)$ is finite.
\end{lemma}

{\bf Proof.} The proof is easy and we left it to the reader.
$\Box$

\begin{lemma} Let $\fa$ be a non-zero ideal of $R$ and $M$ a
$ZD$-module. If $M$ is $\fa$-torsion free, then $\fa$ contains a
nonzero divisor on $M$.
\end{lemma}

{\bf Proof.} Since $M$ is $ZD$-module,  there are prime ideals
$\fp_1,\fp_2,\dots ,\fp_n$ in $\Ass_RM$ such that
$\Z_R(M)=\cup_{i=1}^n\fp_i$. Because, $M$ is $\fa$-torsion free,
it follows that $\fa$ is not contained in any associated prime
ideal of $M$. Thus, by Prime Avoidance Theorem, $\fa$ is not
contained in $\Z_R(M)$. $\Box$

For an $R$-module $M$, the Goldie dimension of $M$ is defined as
the cardinal of the set of indecomposable submodules of $\E(M)$,
which appear in a decomposition of $\E(M)$ into direct sum of
indecomposable submodules. We shall use $\Gdim M$ to denote the
Goldie dimension of $M$. For a prime ideal $\fp$, let
$\mu^0(\fp,M)$ denote the 0-th Bass number of $M$ with respect to
prime ideal $\fp$. It is known that $\mu^0(\fp,M)>0$ if and only
if $\fp\in \Ass_RM$. It is clear by the definition of Goldie
dimension that $\Gdim M=\sum_{\fp\in \Spec R}\mu^0(\fp,M)$. Having
this in mind, we introduce the following generalization of the
notion of Goldie dimension.

\begin{definition} Let $\fa$ be an ideal of $R$. For an $R$-module
$M$, we define $\fa$-relative Goldie dimension of $M$ as
$\Gdim_{\fa}M:=\sum_{\fp\in \V(\fa)}\mu^0(\fp,M)$. Here $\V(\fa)$
denotes the set of prime ideals of $R$ which are containing $\fa$.
\end{definition}

Obviously, if $\fa$ is the zero ideal, then $\Gdim_{\fa}M=\Gdim
M$. Also, it is clear that the Goldie dimension of any Noetherian
module as well as any Artinian module is finite.

\begin{lemma} Let $\fa$ be an ideal of $R$ and $M$ an
$R$-module. Then $\Gdim_{\fa}M=\Gdim \Gamma_{\fa}(M)$.
\end{lemma}

{\bf Proof.} Let $\fp$ be a prime ideal of $R$. By [{\bf 7},
Theorem 18.4], each element of $\E(R/\fp)$ is annihilated by some
power of $\fp$ and for each element $r\in R\smallsetminus \fp$,
the multiplication by $r$ induces an automorphism of $\E(R/\fp)$.
Therefore, it follows that $\E(R/\fp)$ is $\fa$-torsion if
$\fa\subseteq \fp$, and $\fa$-torsion free otherwise. Hence
$\Gamma_{\fa}(\E(M))=\oplus_{\fp\in
\V(\fa)}\mu^0(\fp,M)\E(R/\fp)$. It is easy to see that
$\Gamma_{\fa}(\E(M))$ is an essential extension of
$\Gamma_{\fa}(M)$. On the other hand $\Gamma_{\fa}(\E(M))$ is an
injective $R$-module by [{\bf 2}, Proposition 2.1.4]. Hence
$\Gamma_{\fa}(\E(M))\cong \E(\Gamma_{\fa}(M))$. Thus
$$\Gdim_{\fa}M=\sum_{\fp\in \V(\fa)}\mu^0(\fp,M)=\Gdim
\Gamma_{\fa}(M).\Box$$

\begin{lemma} Let $\fa$ be an ideal of $R$ and $M$  a $ZD$-module. The
following are equivalent:\\
i) $\Gdim_{\fa}M$ is finite.\\
ii) $\Gdim_{\fa R_{\fp}} M_{\fp}$ is finite for any prime ideal
$\fp$ of $R$.\\
iii) $\Gdim_{\fa R_{\fp}} M_{\fp}$ is finite for any prime ideal
$\fp$ which is maximal in $\Ass_RM$.\\
\end{lemma}

{\bf Proof.} First we show that (i) implies (ii). Let $\fp$ be a
prime ideal of $R$ and $S$ a multiplicatively closed subset of
$R$. It follows by [{\bf 2}, Lemma 10.1.12], that if $S\cap
\fp=\emptyset$, then the $S^{-1}R$-modules $S^{-1}(\E(R/\fp))$ and
$\E(R/\fp)$ are isomorphic. Also, if $S\cap \fp\neq\emptyset$ we
can easily deduce that $S^{-1}(\E(R/\fp))=0$. Thus, we have
$$S^{-1}(\E(\Gamma_{\fa}(M)))\cong \oplus_{\fp\in \V(\fa), \fp\cap
S=\emptyset}\mu^0(\fp,M)\E(R/\fp).$$ On the other hand for any
$R$-module $N$, it follows by [{\bf 2}, Corollary 11.1.6], that as
an $S^{-1}R$-module $S^{-1}(\E(N))$ is isomorphic to the
$\E_{S^{-1}R}(S^{-1}N)$. Thus
$$\E_{S^{-1}R}(\Gamma_{\fa S^{-1}R}(S^{-1}M))\cong \E_{S^{-1}R}(S^{-1}
(\Gamma_{\fa}(M))\cong \oplus_{\fp\in
\V(\fa), \fp\cap
S=\emptyset}\mu^0(\fp,M)\E_{S^{-1}R}(S^{-1}R/S^{-1}\fp).$$ This
shows that $\Gdim_{\fa S^{-1}R}S^{-1}M\leq \Gdim_{\fa}M$.
Therefore (i) implies (ii), as required.

Clearly, (ii) implies (iii). Next, we prove that (iii) implies
(i). Let $\{\fq_1,\fq_2,\dots ,\fq_n \}$ be the set of all prime
ideals with the property being maximal in $\Ass_RM$. Note that by
Lemma 2.3, this set is finite. Fix $1\leq i\leq n$. As shown in
the proof of Lemma 2.6, we have
$\E(\Gamma_{\fa}(M))=\oplus_{\fp\in
\V(\fa)}\mu^0(\fp,M)\E(R/\fp)$. Thus

$$\begin{array}{llll} \E_{R_{\fq_i}}(\Gamma_{\fa R_{\fq_i}}(M_{\fq_i}))& \cong
(\E(\Gamma_{\fa}(M)))_{\fq_i}\\ &\cong (\oplus_{\fp\in
\V(\fa)}\mu^0(\fp,M)\E(R/\fp))_{\fq_i} \\ &\cong
\oplus_{\fa\subseteq \fp\subseteq
\fq_i}\mu^0(\fp,M)\E_{R_{\fq_i}}(R_{\fq_i}/\fp
R_{\fq_i}).\end{array}$$ Hence we have $\Gdim_{\fa
R_{\fq_i}}M_{\fq_i}=\sum_{\fa\subseteq \fp\subseteq
 \fq_i}\mu^0(\fp,M)$, and so

$$\begin{array}{llll} \Gdim_{\fa R_{\fq_i}}M_{\fq_i}& =\sum_{\fa\subseteq
\fp\subseteq \fq_i}\mu^0(\fp,M)\\ &\leq \Gdim_{\fa}M \\ &\leq
\sum_{i=1}^n(\sum_{\fa\subseteq \fp\subseteq \fq_i}\mu^0(\fp,M))\\
& =\sum_{i=1}^n(\Gdim_{\fa R_{\fq_i}}M_{\fq_i}).\end{array}$$ This
concludes the proof. $\Box$

\section{Artinianness  of local cohomology modules}

In [{\bf 8}, Theorem 1.3], Melkersson proved that an $\fa$-torsion
module $M$ is Artinian if and only if $0:_M\fa$ is Artinian. In
this section, we use this result to deduce several results
concerning Artinianness  of local cohomology of $ZD$-modules.

\begin{theorem} Let $\fa$ be a non-zero ideal of $R$ and $M$ a
$ZD$-module. The following are equivalent:\\
i) $\Gamma_{\fa}(M/N)$ is Artinian for any submodule $N$ of
$M$.\\
ii)$H^i_{\fa}(M/N)$ is Artinian for any submodule $N$ of $M$ and
all $i\geq 0$.
\end{theorem}

{\bf Proof.} It is clear that (ii) implies (i). Next, we show that
(i) implies (ii) by using induction on $i$. The claim for $i=0$
holds by the assumption. Assume that $i>0$ and that the assertion
holds for $i-1$. Thus $H^{i-1}_{\fa}(M/N)$ is Artinian for all
submodules $N$ of $M$. Let $N$ be a submodule of $M$ and $X=M/N$.
Because $H^i_{\fa}(X)\cong H^i_{\fa}(X/\Gamma_{\fa}(X))$, we may
assume that $X$ is $\fa$-torsion free. Note that any quotient of a
$ZD$-module is also a $ZD$-module. Since $X$ is $\fa$-torsion
free, by  Lemma 2.4, it follows that $\fa$ contains an element $r$
which is nonzero divisor on $X$. The exact sequence $$0\lo X
\overset{r}\lo X\lo X/rX\lo 0,$$ induces an exact sequence
$$H^{i-1}_{\fa}(X/rX)\lo H^i_{\fa}(X) \overset{r}\lo
H^i_{\fa}(X).$$ By inductive hypothesis $H^{i-1}_{\fa}(X/rX)$ is
Artinian, so that by using the above exact sequence, we deduce
that $(0:_{H^i_{\fa}(X)}r)$ is Artinian.  Since $H^i_{\fa}(X)$ is
$Rr$-torsion, the conclusion follows by  [{\bf 8}, Theorem 1.3].
$\Box$

When $\dim R/\fa=0$, we may strengthen Theorem 3.1 as follows.

\begin{proposition} Let the situation be as in Theorem 3.1. In
addition assume that $\dim R/\fa=0$. The following are
equivalent:\\
i) $\Gamma_{\fa}(M/N)$ is Artinian for any submodule $N$ of
$M$.\\
ii) $\fa$-relative Goldie dimension of any quotient of $M$ is
finite.\\
iii)$H^i_{\fa}(M/N)$ is Artinian for any submodule $N$ of $M$ and
all
$i\geq 0$.\\
iv) $H^i_{\fa R_{\fp}}(M_{\fp}/N_{\fp})$ is Artinian for any
submodule $N$ of $M$, any prime ideal $\fp$ of $R$ and all $i\geq
0$.
\end{proposition}

{\bf Proof.} In view of Theorem 3.1, Lemma 2.6 and Lemma 2.7, it
suffices to show that an $\fa$-torsion module $M$ is Artinian if
and only if its Goldie dimension is finite. Assume that $M$ is a
$\fa$-torsion module. Then $\Ass_RM\subseteq \V(\fa)$. On the
other hand, because $\dim R/\fa=0$, it turns out that $\V(\fa)$ is
a finite set consisting of maximal ideals. It is clear that if $M$
is Artinian, then the Goldie dimension of $M$ has to be finite.
Conversely, suppose that the Goldie dimension of $M$ is finite.
Then
$$\sum_{\fp\in \Ass_RM}\mu^0(\fp,M)\leq\sum_{\fp\in
\V(\fa)}\mu^0(\fp,M)<\infty.$$ Thus $\E(M)$ is direct sum of
finitely many Artinian modules. $\Box$

Let $\fa$ be an ideal of a local ring $(R,\fm)$ and let $M$ be a
finitely generated $R$-module of dimension $d$. By [{\bf 2},
Theorem 7.1.6], $H^d_{\fa}(M)$ is Artinian. Also, it is known that
if $\dim R/\fa=0$, then $H^i_{\fa}(M)$ is Artinian for all $i\geq
0$. Next, we provide a far reaching generalization of these facts.

\begin{corollary} Let $\fa$ be an ideal of $R$ and $M$ a
$ZD$-module. Assume that $\fa$-relative Goldie dimension of any
quotient of
$M$ is finite. We have the following.\\
i) If $\dim R/\fa=0$, then $H^i_{\fa}(M)$ is Artinian for all
$i\geq 0$.\\
ii) If $d=\dim M$ is finite, then $H^d_{\fa}(M)$ is Artinian.
\end{corollary}

{\bf Proof.} i) is clear by Proposition 3.2.\\
ii) We use induction on $d$. Suppose $d=0$. Then every associated
prime ideal of $M$ is maximal and so $\E(\Gamma_{\fa}(M))$ is a
direct sum of a finitely many $\E(R/\fm)$, where $\fm$'s are
maximal ideals of $R$. Hence $\Gamma_{\fa}(M)$ is Artinian.

Now, we assume that $d>0$ and that the claim holds for $d-1$.
Similar to the proof of Theorem 3.1, we may assume that $M$ is
$\fa$-torsion free. Thus we can choose an element $r\in \fa$,
which is nonzero divisor on $M$. From the exact sequence $$0\lo M
\overset{r}\lo M\lo M/rM\lo 0,$$ we deduce the exact sequence
$$H^{d-1}_{\fa}(M/rM)\lo H^d_{\fa}(M) \overset{r}\lo
H^d_{\fa}(M),$$ of local cohomology modules. Since $r$ is a
nonzero divisor on $M$, we have $\dim M/rM\leq d-1$. Hence, it
follows from inductive hypothesis or Grothendieck's Vanishing
Theorem [{\bf 2}, Theorem 6.1.2] that $H^{d-1}_{\fa}(M/rM)$ is
Artinian. Therefore by using [{\bf 8}, Theorem 1.3], we deduce
that $H^d_{\fa}(M)$ is Artinian. $\Box$

Next, we bring an example to show that there is a non-finitely
generated $ZD$-module $M$, and an ideal $\fa$ of $R$, such that
$\fa$-relative Goldie dimension of any quotient of $M$ is finite.

\begin{example} i) Let $M$ be a Matlis reflexive $R$-module and
$\fa$ an arbitrary ideal of $R$. Then it follows, by Example
2.2(iii), that $M$ is a $ZD$-module and that $\fa$-relative Goldie
dimension of any quotient of $M$ is finite.

ii) Let $\fm,\fn$ be two distinct maximal ideals of a ring $R$.
Put $M=\oplus_{i\in\mathbb{N}}R/\fm$ and $\fa=\fn$. Then $M$ is a
$ZD$-module and $\fa$-relative Goldie dimension of any quotient of
$M$ is finite. Also, note that by Example 2.2(iii), $M$ is not
Matlis reflexive.
\end{example}


\end{document}